\documentclass[12pt]{amsart}
\usepackage{amssymb}

\usepackage{amscd}

\newtheorem{theorem}{Theorem}
\newtheorem{corollary}{Corollary}
\newtheorem{lemma}{Lemma}
\newtheorem{remark}{Remark}
\newtheorem{definition}{Definition}
\newtheorem{example}{Example}

\begin{document}
\author{Raymond T. Hoobler}
\address{City College of New York\\
New York, New York 10031}
\title{The Merkuriev-Suslin theorem for any semi-local ring}
\maketitle

\begin{abstract}
We introduce here a method which uses \'{e}tale neighborhoods to
extend results from smooth semi-local rings to arbitrary semi-local
rings $A$ by passing to the henselization of a smooth presentation
of $A$. The technique is used to show that \'{e}tale cohomology of
$A$ agrees with Galois cohomology, the Merkuriev-Suslin theorem
holds for $A$, and to describe torsion in $K_2(A)$.
\end{abstract}

 We introduce here a method which uses
\'{e}tale neighborhoods to extend results from smooth semi-local
rings to arbitrary semi-local rings. Three applications are given.
In the first and last, $A$ is a connected, semi-local ring
containing a field $k$ while the second application holds for any
connected, semi-local ring $A$.

1) Let $X=Spec\,(A)$, and let $n$ an integer with $(n,char(k))=1$.
If $F$ is a finite, locally constant sheaf of $\mathbb{Z}/n$-modules
for the \'{e}tale site on $X$, we show
\begin{equation*}
H^{p}(G(A_{s}/A),F(A_{s}))\underset{\approx }{\rightarrow }H^{p}(X,F)
\end{equation*}
where $A_{s}$ is the separable closure of $A,$ the left hand side is the
Galois cohomology of $A$ with coefficients in the $G(A_{s}/A)$-module $%
F(A_{s})$ and the right hand side is the \'{e}tale cohomology group of the
semi-local scheme $X$ with coefficients in $F$.

2) We extend the Merkuriev-Suslin theorem to a connected, semi-local
ring $A$; that is, for $n$ relatively prime to the residue
characteristics of $A$, the Galois symbol map
\begin{equation*}
K_{2}(A)/n\underset{\approx }{\rightarrow }H^{2}(A,\mathbb{Z}/n(2))
\end{equation*}
is an isomorphism where, as usual, $\mathbb{Z}/n(i)=\mu
_{n}^{\otimes i}$. Since this implies the cup product map is
surjective, we conclude that any Azumaya algebra of order $n$ in
$Br(A)$ is similar to a tensor product of symbol algebras if $A$
contains a primitive $n^{th}$ root of unity.

3) We extend Suslin's computation of the $\ell $-primary component
of the torsion in $K_{2}(k)$, $k$ a field, to $A$.

Fix notation as follows. $H^{p}(X,F)$ (or, if $X=Spec\,(A)$, $H^{p}(A,F)$)
denotes the \'{e}tale cohomology group of $X$ with coefficients in the
\'{e}tale sheaf $F$. For connected $A$, Galois cohomology will be indicated
by $H^{p}(\pi _{1}(A,x),M)$ where $\pi _{1}(A,x)$ is the algebraic
fundamental group of $A$ with base point $x:Spec\,(k_{s})\rightarrow
Spec\,(A)$ or, equivalently, by $H^{p}(G(A_{s}/A),M)$ ($\pi _{1}(A,x)=\ $$%
G(A_{s}/A)$, the Galois group of the separable closure of $A$, where the
base point is defined by embedding the residue field of $A$ into a separably
closed field) when $M$ is an abelian group equipped with a continuous
action. $X^{\left( d\right) }$ stands for the set of generic points $x$ of
irreducible components of $X$ of codimension $d$. If $M$ is a fixed abelian
group, we let $M/n$ stand for $M/nM$ , $M\{\ell \}$ stand for the $\ell $%
-primary component of $M$, and $Div(M)$ stand for the maximal divisible
subgroup of $M $ .

I would like to thank Srinivas, Bruno Kahn, and Wayne Raskind for
several very helpful discussions on these and related matters and
Chuck Weibel for the last application and unending patience. I also
want to express my delight and appreciation to the Tata Institute of
Fundamental Research for providing such a fine location and so much
stimulation.

\section{Local functors}

Let $\left( A,I\right) $ be a pair consisting of a commutative ring with $1$
and an ideal $I\subset A.$ An \'{e}tale neighborhood of $\left( A,I\right) $
is a pair $(A^{\prime },I^{\prime })$ and an etale map $u:A\rightarrow
A^{\prime }$ such that $u(I)A^{\prime }=I^{\prime }$ and $u$ induces an
isomorphism $\overline{u}:A/I\rightarrow A^{\prime }/I^{\prime }$.
Geometrically etale neighborhoods of a closed set $W\subset X=Spec(A)$ look
like a cartesian diagram where $\pi $ is etale:
\begin{equation*}
\begin{array}{lll}
W^{\prime } & \subset & X^{\prime } \\
\downarrow \thickapprox &  & \downarrow \pi \\
W & \subset & X
\end{array}
\end{equation*}
The set of all \'{e}tale neighborhoods of $I$ in $A$ is a filtered category
which we denote $Et(I)$. If we assume that $I\subseteq rad(A),$ then $Et(I)$
is used to define the \textit{henselization} of the pair $(A,I)$ \cite
{Raynaud} by setting
\begin{equation*}
A_{/I}{}^{h}:=\underset{(A^{\prime },I^{\prime })\in Et(I)}{\underrightarrow{%
\lim }}A^{\prime }.
\end{equation*}
Then the pair $(A_{/I}^{h},I^{h})$ is a hensel pair where $%
I^{h}=IA_{/I}^{h}. $ If $I=\mathfrak{m}$ is maximal, we get the
hensel local ring $A_{\mathfrak{m}}^{h}. $

Let $\mathcal{C}$ be a category containing $Et(I)$ as a full subcategory and
$A_{/I}^{h}$. We introduce the definition of a \textit{local} functor.

\begin{definition}
A covariant functor $F:\mathcal{C}\rightarrow Ab$ is said to be local (for
the \'{e}tale topology) with respect to a closed set $V(I)$ in $Spec(A)$ if


\begin{enumerate}
\item  $I$ is a radical ideal in $A$, $Et(I)$ is a full subcategory of $\mathcal{C}$%
, and $A_{/I}^{h}\in \mathcal{C}$

\item  the natural map
\begin{equation}
\underset{(A^{\prime },I^{\prime })\in Et(I)}{\underrightarrow{\lim }}%
F(A^{\prime })\rightarrow F(A/I)  \label{1}
\end{equation}

is an isomorphism.
\end{enumerate}

\end{definition}

We say that $F$ is local for the \'{e}tale topology on $X$ if, for all
affine open subsets $Spec(A)\subset X$ and all closed sets $V(I)\subset
Spec(A),$ $F$ is local with respect to the closed set $V(I)\subset Spec(A).$
In many cases of importance the limit condition above can be replaced with a
condition involving the henselization of $A$ along $I$. Recall that a
covariant functor $F:\mathcal{C}\rightarrow Ab$ is said to be locally of
finite presentation if for any filtered inductive limit $A=\lim\limits_{%
\overrightarrow{i\in I}}A_{i}$, the natural map
\begin{equation*}
\lim\limits_{\overrightarrow{i\in I}}F(A_{i})\rightarrow F(A)
\end{equation*}
is an isomorphism.

For such functors we are only interested in a \textit{hensel pair condition}.

\begin{definition}
A functor $F:\mathcal{C}\rightarrow Ab$ is said to satisfy the hensel pair
condition for $I$, resp. epic hensel pair condition for $I$, if $%
F(A)\rightarrow F(A/I)$ is an isomorphism, resp. an epimorphism, for any
hensel pair $(A,I).$ The functor satisfies the hensel pair condition, resp.
epic hensel pair condition if it does so for any hensel pair $(A,I)\in
\mathcal{C}.$
\end{definition}

Then if $F$ is locally of finite presentation and $F$ satisfies the hensel
pair condition, $F$ is local for the closed set $V(I)$. Thus in applying
this definition we will first verify that $F$ is locally of finite
presentation and then that $F(A^{h})\rightarrow F(A^{h}/I^{h})$ is an
isomorphism when $(A^{h},I^{h})$ is a hensel pair.

We are primarily interested in three examples of functors satisfying the
hensel pair condition.

\begin{example}\label{example}
\quad

\begin{enumerate}
\item  $H^{i}(-,F)$ where $F$ is a locally constant sheaf of $\mathbb{Z}/n$
modules on $Spec(A)_{et}$ (\cite{Strano} or \cite{gabberhensel})\label{ex1}

\item  $K_{i}(-;\mathbb{Z}/n)$ where $n$ is relatively prime to the residue
characteristics of $A$ (\cite{gabber})\label{ex2} and

\item  $K_{2}(-)/n$ where $n$ is relatively prime to the residue
characteristics of $A.\label{ex3}$
\end{enumerate}
\end{example}

The last example is easily seen to satisfy the hensel pair condition.
Surjectivity of $K_{2}(A)/n\rightarrow K_{2}(A/I)/n$ for a hensel pair $(A,I)
$ follows immediately since $K_{2}$ is generated by symbols, and injectivity
follows by a straightforward calculation done in \cite[Appendix]{Kahn} or
from Gabber's theorem \cite{gabber}. Alternatively we could use the
universal coefficient theorem, the second example, and note that the third
term in the universal coefficient theorem is the $n-$torsion in $K_{1}$
which is the $n^{th}$ roots of unity. These examples will be discussed at
greater length in the applications.

Our applications are a result of the following theorem.

\begin{theorem}
\label{comparision}\label{fund}Let $A$ be a semi-local ring essentially of
finite type over a base ring $k$. Suppose $A\cong B/J$ where $B$ is a smooth
semi-local ring essentially of finite type over $k$. Let $C$ be a category
of semi-local ring extensions of $B$ containing $A$, $B_{/J}^{h}$ , and the
full subcategory $Et(J)$. Suppose $F_{1},F_{2}:C\rightarrow Ab$ are two
covariant functors which are either both local for $V(J)$ or $F_{2}$ is
local for $V(J)$ and $F_{1}$ is locally of finite presentation and satisfies
the epic hensel pair condition. Let $\phi :F_{1}\rightarrow F_{2}$ be a
natural transformation such that $\phi (B^{\prime })$ is an isomorphism if $%
B^{\prime }$ is smooth and essentially of finite type over $k$. Then $\phi
(A)$ is an isomorphism.
\end{theorem}

\begin{proof}
If $A$ is essentially of finite type over $k$, then $A$ has a presentation
as $B/J$ where $B$ is essentially of finite type and smooth over $k$. Thus $%
\phi (B)$ is an isomorphism as is $\phi (B^{\prime })$ for any $B^{\prime
}\in Et(J)$ and so, by (\ref{1}), $\phi (A)$ must be an isomorphism in the
first case. In the second case, $F_{1}(B_{/J}^{h})\rightarrow F_{1}(A)$ must
be an isomorphism since $F_{2}$ is local for $V(J)$ and $\phi (B_{/J}^{h})$
is an isomorphism.
\end{proof}

\begin{corollary}
In the above situation $F_{1}$ is local for $V(J)$ if it is locally of
finite presentation and satisfies the epic hensel pair condition.
\end{corollary}

\section{Applications}

We have three applications of this perspective. They all rely on first
establishing an isomorphism for semi-local rings smooth over a base scheme
and then using the examples above to extend the result to arbitrary
semi-local schemes over the base.

\subsection{Galois cohomology = \'{E}tale cohomology}

As a first application we consider the relationship between
\'{e}tale and Galois cohomology of semi-local rings. While it is a
generally accepted fact for smooth local rings over a field, the
details do not, as far as I know, appear in the literature although
there is an argument, due to Bruno Kahn, when the generalized Kato
conjecture holds \cite{Kahn}. We begin by reviewing and
reinterpreting Galois cohomology and then proving the isomorphism in
this case.

Let $S$ be an arbitrary connected scheme. Define a site $S_{ret}$ by
considering the category of schemes $T\rightarrow S$ which are finite and
\'{e}tale (= rev\^{e}tement \'{e}tale) over $S$. A covering morphism is a
surjection of schemes over $S$, $T_{2}\rightarrow T_{1}$, which will
necessarily be finite and \'{e}tale. This is the same as considering the
class $E=(ret)$ of all finite, \'{e}tale morphisms
\cite[Chapter II, Section 1]{Milne}. $S_{ret}$ is then the small site $%
(E/S)_{E}$ where $E$ consists of ''rev\^{e}tements \'{e}tales'' as
above. This site is discussed in \cite[Chapter 1]{Milne} where it is
called $FEt$ and used to describe the fundamental group of $S$.
Cohomology in $S_{ret}$ can be identified with Galois cohomology
according to the following recipe.

Suppose $F$ is a sheaf on $S_{ret}$. Then $F$ takes finite, disjoint unions
of schemes to direct products. Thus the discussion in
\cite[Chapter III, Example 2.6]{Milne} applies to the covering $T\rightarrow
S$ in $S_{ret}$ where $T$ is Galois over $S$ with group $G$. Hence
\begin{equation}  \label{Cech}
\check H{}^p(T/S,F)\cong H^p(G,F(T))
\end{equation}
where $F(T)$ is a left $G$-module via the action of $G$ on $T$.

Let $\iota :\mathbf{S}(S_{ret})\rightarrow \mathbf{P}(S_{ret})$ be the
forgetful functor that regards a sheaf as a contravariant functor defined on
$S_{ret}$. When we pass to the limit over all coverings $T/S$ in $S_{ret}$ ,
the left hand side of (\ref{Cech}) becomes $\check{H}{}^{p}(S_{ret},\iota F)$%
. If we fix a geometric point $x:Spec(k_{s})\rightarrow S$ where $k_{s}$ is
a separably closed field, then the theory of the fundamental group expresses
any finite, \'{e}tale covering $T\rightarrow S$ as a quotient of a Galois
covering. Consequently the right hand side of (\ref{Cech}) is $%
\underrightarrow{\lim }\;H^{p}(G_{i},F(T_{i}))$ where $G_{i}$ ranges through
the finite quotients of $\pi _{1}(S,x)$ and $T_{i}\rightarrow S$ is the
corresponding \'{e}tale covering with group $G_{i}$. Let $F(S_{s})$ denote
the abelian group $\underrightarrow{\lim }\;F(T_{i})$. If $F$ is a locally
constant sheaf on $S$, then the subgroup of $\pi _{1}(S,x)$ defining the
covering $T\rightarrow S$ such that $F\mid _{T}$ is constant is of finite
index in $\pi _{1}(S,x)$ and acts trivially on $F(S_{s})$. Hence $F(S_{s})$
has a continuous $\pi _{1}(S,x)$-action. In general this need not be the
case, but we have enough to conclude that for $F$ a locally constant sheaf
on $S_{ret}$,
\begin{equation}
H^{p}(S_{ret},F)\cong H^{p}(\pi _{1}(S,x),F(S_{s}))  \label{ciso}
\end{equation}
(Note that a finite, locally constant sheaf for the \'{e}tale topology is
the same as a finite, locally constant sheaf on $S_{ret}$.) In view of the
isomorphism (\ref{Cech}) on $S_{ret}$, this is essentially the statement
that sheaf cohomology coincides with \v{C}ech cohomology. Corollary 2.5 of
Chapter III \cite{Milne} gives sufficient conditions for this; namely, for
every surjection $F\rightarrow F^{\prime \prime }$ of sheaves, the map
\begin{equation*}
\underrightarrow{\lim }\prod F(U_{i_{0}...i_{p}})\rightarrow
\underrightarrow{\lim }\prod F^{\prime \prime }(U_{i_{0}...i_{p}})
\end{equation*}
is surjective where the limit is over all coverings $\left\{ \coprod
U_{i}\rightarrow S\right\} $ of $S$ and $U_{i_{0}...i_{p}}=U_{i_{0}}\times
_{S}\cdots \times _{S}U_{i_{p}}.$ But if $T\rightarrow S$ is a Galois
covering of $S$ with group $G$, then $F(T_{i_{0}...i_{p}})=F(\coprod_{G^{%
\times p}}T)=\prod_{G^{\times p}}F(T)$. Thus taking limits of $%
F(T_{i_{0}...i_{p}})$ over all coverings means taking limits over all
coverings of products of copies of $F(T)$ and so surjectivity is immediate.

This discussion reduces our first application to the following theorem about
the change of sites morphism $\tau :S_{et}\rightarrow S_{ret}$.

\begin{theorem}
Let $A$ be a connected semi-local ring containing a field $k$. If
$F$ is a finite, locally constant sheaf of $\mathbb{Z}/n$-modules
for the \'{e}tale site on $A$ where $(n,char(k))=1$, then
\begin{equation}
H^{p}(A_{ret},\tau _{*}F)\rightarrow H^{p}(A_{et},F)  \label{iso1}
\end{equation}
is an isomorphism for all p.
\end{theorem}

\begin{proof}
We wish to apply our extension theorem. We begin by treating the case of a
smooth, semi-local ring $R$ essentially of finite type over $k$, an
algebraically closed field first. We will show that $R^{q}\tau _{*}F=0$ for $%
q>0$ and any finite, locally constant sheaf $F$ of
$\mathbb{Z}/n$-modules. Suppose $x\in H^{q}(R,F)$ with $q>0.$ We may
assume $F=C_{X}$ is a constant sheaf of $\mathbb{Z}/n$-modules with
value $C$ since $F$ becomes constant after a finite, etale extension
of $R$. By assumption there is a smooth, connected
variety $X$ over $k$, closed points $t_{i}\in X,\;1\leq i\leq m$, such that $%
R\cong \mathit{O}_{X,t_{1}...t_{m}}$ and an element $x^{\prime }\in
H^{q}(X,C_{X})$ such that $x^{\prime }\mid _{Spec(\mathit{O}%
_{X,t_{1}...t_{m}})}=x$. We construct a finite, \'{e}tale covering $\rho
:X^{\prime }\rightarrow X$ such that $\rho ^{*}(x^{\prime })=0$ using
induction on $\dim (X)$ and the existence of Artin neighborhoods. If $\dim
(X)=0$, the assertion is obvious. If $\dim (X)=d$, the existence of an Artin
neighborhood relative to $k$ means there is a diagram (\cite[p. 117]{Milne}
or \cite{sga4})
\begin{equation*}
\begin{array}{ccccc}
X & \;\overset{j}{\rightarrow } & \overline{X} & \overset{i}{\leftarrow } & Y
\\
& f\searrow & \downarrow \overline{f} & \swarrow g &  \\
&  & X_{0} &  &
\end{array}
\end{equation*}
in which

\begin{enumerate}
\item  $j$ is an open immersion, dense in each fibre and $Y=\overline{X}-X;$

\item  $\overline{f}$ is smooth and projective with geometrically
irreducible fibres of dimension one;

\item  $g$ is finite and \'{e}tale and each fibre of $g$ is non-empty.
\end{enumerate}

(Although the statement of the existence of Artin neighborhoods given in
\cite{Milne} refers only to a local ring, the proof as given in
\cite[Expose XI, Section 3]{sga4} clearly extends to the semi-local case. In
fact, the proof of Proposition 3.3, the key statement, is footnoted to that
effect. Since base extension of a local ring from an arbitrary field $k$ to
its separable closure $k_{s}$ may produce a semi-local ring we need this
more general result.)

Now $(Y,\overline{X})$ is a smooth $X_{0}-$pair. The Gysin sequence
\cite[Chapter VI, Corollary 5.3 and Remark 5.4 (a)]{Milne} becomes $%
\overline{f}_{*}C_{\overline{X}}\cong f_{*}C_{X},$%
\begin{equation*}
0\rightarrow R^{1}\overline{f}_{*}C_{\overline{X}}\rightarrow
R^{1}f_{*}C_{X}\rightarrow g_{*}(C_{Y}\otimes T_{Y/\overline{X}})\rightarrow
R^{2}\overline{f}_{*}C_{\overline{X}}\rightarrow R^{2}f_{*}C_{X}\rightarrow
0,
\end{equation*}
and $R^{p}\overline{f}_{*}C_{\overline{X}}\cong R^{p}f_{*}C_{X}$ for $p>2$.
The base change maps for the proper morphisms $\overline{f}$ and $g$ are
isomorphisms. Consequently the base change map for $f$ is also an
isomorphism. Since the fibres of $f$ are non-complete curves, we find $%
R^{p}f_{*}C=0$ for $p>1$. Moreover $\overline{f}$ is smooth and so $R^{1}%
\overline{f}_{*}C_{\overline{X}}$ is finite and locally constant on $X_{0}$
as is $g_{*}(C_{Y}\otimes T_{Y/X})$ since $g$ is an etale covering space and
$T_{Y/X}$ is locally constant on $Y$. Therefore we conclude that $%
R^{p}f_{*}C_{X}$ is a finite, locally constant sheaf of
$\mathbb{Z}/n-$modules for $p=0$ or $1$.

This information shows that the Leray spectral sequence for $f$ degenerates
into the long exact sequence
\begin{equation}
\begin{split}
 ...\rightarrow H^{q}(X_{0},f_{*}C_{X})\rightarrow
H^{q}(X,C_{X})\rightarrow & H^{q-1}(X_{0},R^{1}f_{*}C_{X})\\
&\rightarrow H^{q+1}(X_{0},f_{*}C_{X})\rightarrow ...
\end{split}
\end{equation}
and we may apply the induction hypothesis. Then for any
\begin{equation*}
y\in H^{q-1}(X_{0},R^{1}f_{*}C_{X}),
\end{equation*}
 there is,
after shrinking $X_{0}$ if necessary, a finite, \'{e}tale covering
space $\rho ^{\prime }:X_{0}^{\prime }\rightarrow X_{0}$ such that
$\rho ^{\prime *}(y)=0$. Now by using the pullbacks of appropriate
coverings of $X_{0}$ to $X$ and a diagram chase, it is a
straightforward matter to produce a finite, \'{e}tale covering $\rho
:X^{\prime }\rightarrow X$ such that $\rho ^{*}(x^{\prime })=0$.

This argument is, of course, the essence of Remark 3.16 \cite[Chapter III]
{Milne} and shows that $R^{q}\tau _{*}F=0$ for $q>0$. Hence the Leray
spectral sequence for the change of site morphism $\tau :X_{et}\rightarrow
X_{ret}$ collapses, and we conclude
\begin{equation*}
H^{p}(R_{ret},\tau _{*}F)\cong H^{p}(R_{\acute{e}t},F)
\end{equation*}

Now suppose that $R$ is a semi-local ring which is smooth and essentially of
finite type over $k$, a not necessarily algebraically closed field. Then for
any scheme $X$ there is a Hochschild-Serre spectral sequence
\cite[Chapter III, Remark 2.21(b)]{Milne}
\begin{equation} \label{HS1}
H^{p}(G(k_{s}/k),H^{q}(X_{s},F))\Rightarrow H^{n}(X,F)
\end{equation}
where $X_{s}=X\times _{k}Spec(k_{s})$. Similarly the theory of algebraic
fundamental groups constructs a profinite group extension
\begin{equation*}
1\rightarrow G((R\otimes _{k}k_{s})_{s}/R\otimes _{k}k_{s})\rightarrow
G(R_{s}/R)\rightarrow G(k_{s}/k)\rightarrow 1
\end{equation*}
and so there is a corresponding Hochschild-Serre spectral sequence
\begin{equation}  \label{HS2}
\begin{split}
H^{p}(G(k_{s}/k),H^{q}(G((R\otimes _{k}k_{s})_{s}/R\otimes
_{k}k_{s}),F((R\otimes _{k}k_{s})_{s})))&\Rightarrow\\
H^{n}(G(R_{s}/R),F(R_{s}))
\end{split}
\end{equation}
The change of site morphism defines a homomorphism $\Psi $ from (\ref{HS1})
to (\ref{HS2}).

Since $F$ is a sheaf of $\mathbb{Z}/n$-modules where $n$ is
relatively prime to char$(k)$, we conclude, as usual, that
$H^{p}(X\times
_{k}Spec(k_{s}),F)\cong H^{p}(X\times _{k}Spec(\overline{k}),F)$ where $%
\overline{k}$ is the algebraic closure of $k$. Consequently, by the first
case, $\Psi $ is an isomorphism on the $E_{2}^{p,q}$ terms. This shows that (%
\ref{iso1}) is an isomorphism when $A$ is essentially smooth over a not
necessarily algebraically closed field $k.$

Some preparation is required in order to apply Theorem
\ref{comparision} to an arbitrary semi-local ring $A$ whether we
work with the small or large \'{e}tale site. Our hypothesis on $F$
shows that $F$ is representable in either the small or large site by
a semi-local scheme $\mathbb{F}$ which is finite and \'{e}tale over
$A.$ Now $A$ is a filtered limit of connected semi-local rings
$A_{i}$ essentially of finite type over $k.$ Consequently there is a
finite, etale covering space $\mathbb{F}_{i_{0}}$ over $A_{i_{0}}$
for a sufficiently large $A_{i_{0}}$ such that
$\mathbb{F}_{i_{0}}\times _{A_{i_{0}}}A\approxeq \mathbb{F}$. We use
$\mathbb{F}_{i_{0}}$ to define a locally constant sheaf $F_{i_{0}}$
over $A_{i_{0}}$ whose restriction to $A$ is $F.$ Next find a
presentation of $A_{i_{0}}$ by a smooth connected
semi-local ring $B_{i_{0}}$ essentially of finite type over $k$ so that $%
A_{i_{0}}\cong B_{i_{0}}/I$. Since any \'{e}tale covering space of $%
A_{i_{0}} $ may be lifted to an \'{e}tale covering space of $B_{i_{0}}$ (by
the theorem of the primitive element all one has to do is lift a separable
polynomial), we may assume there is a finite, locally constant sheaf $G$ on $%
Spec(B_{i_{0}})$ whose restriction to $Spec(A)$ is $F$.

Let $C$ be the category of connected semi-local rings over $B$, $%
F_{1}=H^{p}(-_{ret},\tau _{*}G)$, and $F_{2}=H^{p}(-,G)$ where now we
consider $G$ on the large site over $Spec(B)$. Note that $G$ is represented
by a finite \'{e}tale covering space $B^{\prime }$ of $B$ in $B_{et}$ and so
$G$ is also represented by $B^{\prime }$ on the big \'{e}tale site over $B$.
Since $B^{\prime }$ also defines a covering in $B_{ret}$ the same
observation applies to the sheaf $\tau _{*}G$ in the big site $B_{ret}$.
After showing that the hypotheses of our extension theorem are met we can
conclude that (\ref{iso1}) is an isomorphism for any $B-$algebra $A$ with
the property that $G\mid _{A}=F$ and $\tau _{*}G\mid _{A}=\tau _{*}F$.

Both $F_{1}$ and $F_{2}$ are locally of finite presentation. For $F_{2}$
this is \cite[Chapter III, Lemma 1.16]{Milne} suitably interpreted. For $%
F_{1}$ we need a different argument. Suppose $R=\underrightarrow{\lim }%
\,R_{i}$ with $R_{i}$ a connected semi-local ring over $B$ where $i\in
\mathit{I}$, a filtered set. Any \'{e}tale covering space of $%
\underrightarrow{\lim }\,R_{i}$ comes from an \'{e}tale covering space
defined over one of the $R_{i}$ and homomorphisms between any two \'{e}tale
covering spaces of $\underrightarrow{\lim }\,R_{i}$ descend to a
homomorphism between covering spaces over one of the $R_{i}$. Thus we have
an isomorphism
\begin{equation*}
\lim_{R_{i}}\lim_{T_{i}/R_{i}}\check{H}{}^{p}(T_{i}/R_{i},G)\rightarrow
\lim_{T/R}\check{H}{}^{p}(T/R,G)
\end{equation*}
since the \v{C}ech cohomology groups depend only on the Galois group of the
covering by (\ref{Cech}) and $G$ is locally constant. Thus, using the Galois
cohomology interpretation of $F_{1}$, we also have
\begin{equation*}
H^{p}(R_{ret},\tau _{*}G)\cong \underrightarrow{\lim }\,H^{p}(R_{i\,ret},%
\tau _{*}G)
\end{equation*}
and so $F_{1}$ is locally of finite presentation.

Now we turn to the hensel pair condition. Suppose $(R,J)$ is a
hensel couple and $G$ is a locally constant sheaf of
$\mathbb{Z}/n-$modules on $Spec(R)$. Then $G$ is represented by a
finite, \'{e}tale group scheme over $R$, and, by Example \ref{ex1},
we have an isomorphism for a hensel couple $(R,J)$ and such a $G$,
\begin{equation*}
H^{p}(R_{\acute{e}t},G)\cong H^{p}(R/J_{\acute{e}t},G\mid _{R/J}).
\end{equation*}
Thus $F_{2}$ is a local functor for the etale topology satisfying the hensel
pair condition. As for $F_{1}$, there is a one-to-one correspondence between
\'{e}tale covering spaces of $R/J$ and \'{e}tale covering spaces of $R$, and
so we have an isomorphism
\begin{equation*}
H^{p}(R_{ret},\tau _{*}G)\cong H^{p}(R/J_{ret},\tau _{*}(G\mid _{R/J})).
\end{equation*}
Hence $F_{1}$ is a local functor for the etale topology and also satisfies
the hensel pair condition. We can now apply Theorem \ref{comparision} to
conclude the proof.
\end{proof}

\begin{corollary}
Let $A$ be a semi-local ring containing a field $k$. Suppose $F$ is
a finite, locally constant sheaf of $\mathbb{Z}/n$-modules where
$(n,char(k))=1$.
Then for any $x\in H^{p}(A,F)$, $p>0$, there is a Galois extension $%
A^{\prime }/A$ such that $x\mid _{A^{\prime }}=0\in H^{p}(A^{\prime },F)$.
\end{corollary}

\subsection{Merkurjev-Suslin theorem}

Our second application is an extension of the Mercuriev-Suslin theorem to an
arbitrary semi-local ring.

For a semi-local ring $A$ and any $n$ which is relatively prime to
the residue characteristics of $A$, Kummer theory provides a natural
isomorphism $K_{1}(A)/n\rightarrow H^{1}(A,\mathbb{Z}/n(1)).$ The
Galois symbol map \cite{Tate}
is obtained by extending this map multiplicatively to $%
K_{i}(A)/n\rightarrow H^{i}(A,\mathbb{Z}/n(i)).$

\begin{theorem}
Let $A$ be a semi-local ring. Then if $n$ is relatively prime to the residue
characteristics of $A$, the Galois symbol map
\begin{equation}
K_{2}(A)/n\rightarrow H^{2}(A,\mathbb{Z}/n(2))
\end{equation}
is an isomorphism.
\end{theorem}

\begin{corollary}
Let $A$ be a semi-local ring, and suppose $n$ is an integer which is
relatively prime to the residue characteristics of $A$. Then the cup product
map
\begin{equation*}
H^{1}(A,\mathbb{Z}/n(1))\times H^{1}(A,\mathbb{Z}/n(1))\rightarrow H^{2}(A,\mathbb{Z}%
/n(2))
\end{equation*}
is onto. In particular, if $A$ contains a primitive $n^{th}$ root of unity,
the $n$-torsion in the Brauer group of $A$ is generated by symbol algebras $%
(a,b)_{n}$.
\end{corollary}

\begin{corollary}
Let $X$ be a scheme, and suppose $n$ is an integer which is relatively prime
to the residue characteristics of $X$. The symbol map defines an isomorphism
of Zariski sheaves
\begin{equation*}
\underline{K_{2}/n}\rightarrow \underline{H^{2}(\mathbb{Z}/n(2)}).
\end{equation*}
\end{corollary}

The proof requires some preliminary material on K-theory with coefficients,
Chern classes, and a discussion of results of Gillet which will be used to
prove the theorem for semi-local rings smooth and essentially of finite type
over a semi-local Dedekind ring.

Fix a ring $A$, a prime $\ell $ that is a unit in $A$, and an
integer $n$ which is at least $2$ if $\ell =2$ (See
\cite[Proposition 2.4]{LLouise} for a detailed discussion of the
case $\ell =2$.) Algebraic K-theory with coefficients,
$K_{r}(A;\mathbb{Z}/\ell ^{n})$, is a bifunctor in $A$ and the
cyclic group $\mathbb{Z}/\ell ^{n}$, and, for a pair of integers
$1\leq i,0\leq j\leq 2i$, there is a theory of Chern classes given
by natural transformations
\begin{equation*}
c_{i,j}^{\ell ^{n}}:K_{2i-j}(A;\mathbb{Z}/\ell ^{n})\rightarrow H^{j}(A,\mathbb{Z}%
/\ell ^{n}(i))
\end{equation*}
satisfying certain properties. Among the properties these objects satisfy
are:

\begin{enumerate}
\item  (universal coefficient theorem) There is a natural short exact
sequence
\begin{equation}
0\rightarrow K_{r}(A)/\ell ^{n}K_{r}(A)\rightarrow
K_{r}(A;\mathbb{Z}/\ell ^{n})\rightarrow \ _{\ell
^{n}}K_{r-1}(A)\rightarrow 0  \label{uct}
\end{equation}

\item  (functoriality) For any $n>m$, the diagram
\begin{equation*}
\begin{array}{ccc}
K_{2i-j}(A;\mathbb{Z}/\ell ^{m}) & \rightarrow &
K_{2i-j}(A;\mathbb{Z}/\ell ^{n})
\\
\downarrow c_{ij}^{\ell ^{m}} &  & \downarrow c_{ij}^{\ell ^{n}} \\
H^{j}(A,\mathbb{Z}/\ell ^{m}(i)) & \rightarrow &
H^{j}(A,\mathbb{Z}/\ell ^{n}(i))
\end{array}
\end{equation*}
defined from the inclusion $\mathbb{Z}/\ell ^{m}\hookrightarrow
\mathbb{Z}/\ell ^{n}$ commutes.

\item  (naturality) For any valuation ring $\mathcal{O}_{v}$ and prime $\ell
$ distinct from the residue characteristic of $\mathcal{O}_{v}$ and any $n$,
we have a commutative diagram
\begin{equation*}
\begin{array}{ccc}
K_{2i-j}(F;\mathbb{Z}/\ell ^{m}) & \overset{\overline{\partial
}}{\rightarrow }
& K_{2i-j-1}(\kappa (v);\mathbb{Z}/\ell ^{n}) \\
\downarrow c_{ij}^{\ell ^{n}} &  & \downarrow (1-i)c_{i-1j-1}^{\ell ^{n}} \\
H^{j}(F,\mathbb{Z}/\ell ^{n}(i)) & \overset{\partial }{\rightarrow }
& H^{j-1}(\kappa (v),\mathbb{Z}/\ell ^{n}(i-1))
\end{array}
\end{equation*}
where $\kappa (v)$ is the residue field of $\mathcal{O}_{v}$ and $F$
is its field of fractions. In particular the Chern classes fit
together to define a map from the Quillen-Gersten complex for
$K_{*}(-;\mathbb{Z}/\ell ^{n})$ to the Bloch-Ogus complex for
$H^{*}(-,\mathbb{Z}/\ell ^{n}(\cdot ))$. (Quillen's
argument in \cite{Quillen} applies unchanged to $K_{i}(-;\mathbb{Z}/m)$ where $%
m=\ell ^{n}$.)
\end{enumerate}

Details of the above properties can be found in \cite{soule}.

Now the other tool we need is a reformulation and extension to a
semi-local ring of a result of Gillet-Levine\cite{gilletk} and
Gillet\cite{gillet}.

\begin{theorem}[Gillet-Levine, Gillet]
Let $B$ be a connected, semi-local ring with quotient field $K$ which is
smooth and essentially of finite type over a semi-local Dedekind ring $D$.
Then, for any integer $n$ relatively prime to the residue characteristics of
$B$ and any $q\geq 0,$ we have a Gersten-Quillen resolution of $K_{q}(B;\mathbb{%
Z}/n):$%
\begin{equation*}
0\rightarrow K_{q}(B;\mathbb{Z}/n)\rightarrow
K_{q}(K;\mathbb{Z}/n)\rightarrow \coprod_{x\in
B^{(1)}}K_{q-1}(\kappa (x);\mathbb{Z}/n)\rightarrow \cdots
\end{equation*}
\end{theorem}

\begin{proof}
The argument of Gillet and Levine and Gillet immediately extends to
semi-local rings as is clear from going through their arguments.
\end{proof}

We also need a weak version of this result for etale cohomology.
This has been demonstrated by Gillet but remains
unpublished\cite{gilletbo}. A proof for $n$ an odd integer
relatively prime to the residue characteristics of $B$ is given in
an appendix at the end of this paper.

\begin{lemma}
\label{Lemma}\label{keylemma}Let $B$ be a connected, semi-local ring which
is smooth and essentially of finite type over a semi-local Dedekind ring, $%
\mathcal{D}$. Let $K$ be the quotient field of $B.$ Then, for any integer $n$
which is relatively prime to the residue characteristic of $\mathcal{D},$
\begin{equation*}
H^{2}(B,\mathbb{Z}/n(2))\rightarrow H^{2}(K,\mathbb{Z}/n(2))
\end{equation*}
is a monomorphism.
\end{lemma}

We will need the universal coefficient theorem to produce a Gersten-Quillen
sequence for $K_{2}(-)/n$ using work of Gillet and Gillet-Levine and then
the compatibility of the Chern class map with the Quillen-Gersten complex
and the Bloch-Ogus complex will do the rest in the smooth case.

\begin{proof}
We deal with a semi-local algebra $A$ over a mixed characteristic Dedekind
ring $\mathcal{D}$ first. We may assume $A$ is connected. Both $K_{2}(-)/n$
and $H^{2}(-,\mathbb{Z}/n(2))$ are locally of finite presentation over $Spec(%
\mathcal{D}),$ and so we may assume $A$ is essentially of finite type over $%
\mathcal{D}.$ Now realize $A\approxeq B/J$ where $B$ is a smooth, semi-local
ring which is essentially of finite type over $\mathcal{O}_{v}.$ Let $%
\mathcal{C}$ be the category of connected semi-local rings essentially of
finite presentation over $\mathcal{D}.$ We wish to apply Theorem \ref{fund}
to this situation with $k=\mathcal{D},$ $F_{1}=K_{2}(-)/n,$ $F_{2}=H^{2}(-,%
\mathbb{Z}/n(2)),$ and $\phi :F_{1}\rightarrow F_{2}$ being the
Galois symbol map.

First observe that $K_{2}(-)/n$ satisfies the epic hensel pair
condition since it is generated by symbols and
$H^{2}(-,\mathbb{Z}/n(2))$ satisfies the hensel pair condition by
Gabber or Strano's result (\cite{gabberhensel} or \cite{Strano}). It
remains to show that $\phi (B)$ is an isomorphism if $B$ is smooth
and essentially of finite type over $\mathcal{D}.$ The exactness of
\begin{equation}\label{relkseq}
K_{2}(B;\mathbb{Z}/n)\hookrightarrow
K_{2}(K;\mathbb{Z}/n)\rightarrow \coprod_{x\in
(SpecB)^{(1)}}K_{1}(\kappa (x);\mathbb{Z}/n)
\end{equation}
where $K$ is the quotient field of $B$ was shown by Gillet
\cite{gillet} using work of Gillet-Levine \cite{gilletk}(see Theorem
2.6).
We see that $K_{1}(\kappa (x);\mathbb{Z}/n)\approxeq K_{1}(\kappa (x))/n,$ $%
K_{2}(B;\mathbb{Z}/n)$ is an extension of $\mu _{n}(B)(=_{n}K_{1}(B))$ by $%
K_{2}(B)/n,$ and similarly for $K_{2}(K;\mathbb{Z}/n).$ We can now
reinterpret (\ref{relkseq}) as the exact sequence of the first line
below
\begin{equation*}
\begin{array}{cccc}
& K_{2}(B)/n\hookrightarrow & K_{2}(K)/n\rightarrow & \coprod_{x\in
(SpecB)^{\left( 1\right) }}K_{1}(\kappa (x))/n \\
& \downarrow & \downarrow \approxeq & \downarrow \approxeq \\
& H^{2}(B,\mathbb{Z}/n(2))\hookrightarrow &
H^{2}(K,\mathbb{Z}/n(2))\rightarrow & \coprod_{x\in (SpecB)^{\left(
1\right) }}H^{1}(\kappa (x),\mathbb{Z}/n(1))
\end{array}
.
\end{equation*}
Here the bottom row is not necessarily exact but is a complex, and the first
map is a monomorphism by Lemma \ref{keylemma}, the first vertical
isomorphism is the Merkuriev-Suslin theorem for fields, and the second
vertical isomorphism is the observation that both groups are isomorphic to $%
\coprod_{x\in (SpecB)^{\left( 1\right) }}\kappa (x)^{*}/\kappa (x)^{*n}.$
The diagram commutes by the naturality condition above. Since the bottom row
is a complex we conclude that the first vertical map is an isomorphism as
desired. Theorem 1.1 now finishes this case.

The case of a semi-local ring containing a field is similar but simpler
since we can use Grayson's version of the Gersten-Quillen sequence \cite
{Grayson} and the Bloch-Ogus sequence \cite{CHK} in place of the argument
involving algebraic K-theory with coefficients and Lemma \ref{keylemma}.
\end{proof}

Suppose the Bloch-Kato conjecture holds for fields; that is, the Galois
symbol map
\begin{equation}
K_{i}^{M}(K)/n\rightarrow H^{i}(K,\mathbb{Z}/n(i))  \label{Kato}
\end{equation}
is an isomorphism where $K_{i}^{M}(K)$ is the $i^{th}$ Milnor
K-group of the field $K$ and $n$ is relatively prime to $char(K)$.
Then the same argument applies since the Gersten-Quillen and
Bloch-Ogus sequences hold for $i$ as well as 2, at least up to
$(i-1)!$ torsion and Gillet's unpublished result. Thus the
Bloch-Kato conjecture for fields in degree $i$ and $i-1$ would show
that (\ref{Kato}) is an isomorphism for smooth local rings
essentially of finite type over $k$. The rest of the argument is
then identical, and so we would conclude that $K_{i}^{M}(A)/n$
$\rightarrow H^{i}(A,\mathbb{Z}/n(i))$ is
an isomorphism if $A$ is a semi-local ring containing $k$ and $%
(n,(i-1)!char(k))=1$.

\subsection{Torsion in K$_{2}$}

The final application was suggested by C. Weibel. Recall that Suslin in \cite
{suslin} used Chern classes to construct an isomorphism
\begin{equation*}
K_{2}(F)\{\ell \}\underset{\approx }{\rightarrow }H^{1}(F,\mathbb{Q}\mathbf{%
_{\ell }/}\mathbb{Z}_{\ell }(2))/Div(H^{1}(F,\mathbb{Q}\mathbf{_{\ell }/}\mathbb{Z}%
_{\ell }(2)))
\end{equation*}
where $F$ is either a field of positive characteristic and $\ell $
is a prime distinct from $char(F),$ or $F$ is a field of finite type
over $\mathbb{Q} $ . Note that Suslin is using continuous Galois
cohomology since he uses the identification
\begin{equation*}
H^{1}(F,\mathbb{Q}\mathbf{_{\ell }/}\mathbb{Z}_{\ell }(2))/Div(H^{1}(F,\mathbb{Q}%
\mathbf{_{\ell }/}\mathbb{Z}_{\ell
}(2)))=H_{cont}^{2}(F,\mathbb{Z}_{\ell }(2))\{\ell \}
\end{equation*}
where continuous cochains must be used to correctly identify $H_{cont}^{2}(F,%
\mathbb{Z}_{\ell }(2))\{\ell \}$ as $H^{1}(F,\mathbb{Q}_{\ell
}/\mathbb{Z}_{\ell })(2))$ modulo its maximal divisible subgroup.

\begin{theorem}
\label{sus}Let $A$ be a semi-local ring containing a field $k$ and
suppose that $A$ is either essentially of finite type over
$\mathbb{Q}$ or $char(k)>0$.
Let $\ell $ be a prime distinct from $char(k)$. Then the Chern class $%
c_{2,1} $ induces an isomorphism
\begin{equation*}
K_{2}(A)\{\ell \}\underset{\approx }{\rightarrow }H^{1}(A,\mathbb{Q}\mathbf{%
_{\ell }/}\mathbb{Z}_{\ell }(2))/Div(H^{1}(A,\mathbb{Q}\mathbf{_{\ell }/}\mathbb{Z}%
_{\ell }(2)))\underset{\approx }{\rightarrow
}H_{cont}^{2}(A,\mathbb{Z}_{\ell }(2))\{\ell \}
\end{equation*}
where $H_{cont}^{2}(A,\mathbb{Z}_{\ell })(2))$ is continuous
\'{e}tale cohomology \cite{Jannsen}.
\end{theorem}

\begin{corollary}
Suppose $A$ is a unibranch, e.g. normal, local ring containing a field $k$
with $char(k)=p>0,$ and let $\ell $ be as above. Then
\begin{equation*}
K_{2}(A)\{\ell \}\underset{\approx }{\rightarrow }H^{1}(A,\mathbb{Q}_{\ell }/%
\mathbb{Z}_{\ell }(2))
\end{equation*}
\end{corollary}

The conclusion of Theorem \ref{sus} is independent of the characteristic of $%
k$, but the argument required in the two cases is different. We begin by
recalling the structure of Suslin's argument in \cite{suslin} and the role
of continuous cohomology. He first observes (in the proof of Proposition 3.8$%
)$ that $c_{2,1}$ appears in a commutative, exact diagram

\begin{equation}
\begin{array}{ccccc}
K_{3}(F)\rightarrow & K_{3}(F;\mathbb{Z}/\ell ^{n}) & \overset{c_{2,1}}{%
\rightarrow } & H^{1}(F,\mathbb{Z}/\ell ^{n}(2)) &  \\
\downarrow &  &  & \downarrow &  \\
H_{cont}^{1}(F,\mathbb{Z}_{\ell }(2))\rightarrow &
H_{cont}^{1}(F,\mathbb{Q}_{\ell }(2)) & \rightarrow &
H_{cont}^{1}(F,\mathbb{Q}_{\ell }/\mathbb{Z}_{\ell }(2))\rightarrow
& H_{cont}^{2}(F,\mathbb{Z}_{\ell }(2))\{\ell \}\rightarrow 0
\end{array}
\label{c21}
\end{equation}
where the right hand vertical arrow comes from the coefficient
sequence defined by multiplication by $\ell ^{n}$ on
$\mathbb{Q}_{\ell }/\mathbb{Z}_{\ell }(2)$. Continuous cohomology is
used here to make the bottom sequence exact (the first arrow is
defined on $\ell $-adic sheaves by sending $x\mapsto \ell ^{-n}x$)
and to identify the image of$H_{cont}^{1}(F,\mathbb{Q}_{\ell }(2)) $
with the maximal divisible subgroup of
$H_{cont}^{1}(F,\mathbb{Q}_{\ell
}/\mathbb{Z}_{\ell }(2))$ and the image of $H_{cont}^{1}(F,\mathbb{Q}_{\ell }/\mathbb{%
Z}_{\ell }(2))$ with the $\ell $-primary torsion in $H_{cont}^{2}(F,\mathbb{Z}%
_{\ell }(2))$. Now by exactness and the universal coefficient theorem $%
c_{2,1}$ factors through $K_{2}(F)\{\ell \}$ as
\begin{eqnarray*}
\overline{c}_{2,1} :&K_{2}(F)\{\ell \}\rightarrow H_{cont}^{1}(F,\mathbb{Q}%
_{\ell }/\mathbb{Z}_{\ell }(2))/Image(H_{cont}^{1}(F,\mathbb{Q}_{\ell }(2))) \\
&\cong H_{cont}^{2}(F,\mathbb{Z}_{\ell }(2))\{\ell \}.
\end{eqnarray*}
Suslin has already shown that $H_{cont}^{1}(F,\mathbb{Z}_{\ell
}(2))$ is finite
when $char(F)>0$ \cite[Corollary 2.8]{suslin}, and so $H_{cont}^{1}(F,\mathbb{Q}%
_{\ell }(2))=0$ and
\begin{equation*}
H_{cont}^{1}(F,\mathbb{Q}_{\ell }/\mathbb{Z}_{\ell }(2))\cong H_{cont}^{2}(F,\mathbb{Z%
}_{\ell }(2))\{\ell \}
\end{equation*}
in this case. He then shows that $\overline{c}_{2,1}:K_{2}(F)\{\ell
\}\rightarrow H_{cont}^{1}(F,\mathbb{Q}_{\ell }/\mathbb{Z}_{\ell
}(2))$ is an isomorphism \cite[Theorem 3.9]{suslin} by an induction
argument. In
characteristic zero $H_{cont}^{1}(F,\mathbb{Q}_{\ell }(2))$ is only zero if $%
F_{0}$, the algebraic closure of $\mathbb{Q}$ in $F$, has only real
embeddings in $\mathbb{C}$, but he can reduce the theorem to $F_{0}$
\cite[Proposition 3.3] {suslin} and here the desired result was
proven earlier by Tate. Thus extending this result to a semi-local
ring containing a field which is essentially of finite type over
$\mathbb{Q}$ in characteristic zero requires first defining the
factorization of $c_{2,1}$ through $\overline{c}_{2,1}$ (which will
require continuous \'{e}tale cohomology \cite{Jannsen}) and then
showing that $\overline{c}_{2,1}$ is an isomorphism.

$c_{2,1}$ is defined for arbitrary semi-local rings $A$ (at least if $n>2$
when $\ell =2$) and Suslin's factorization argument applies which allows us
to define
\begin{equation}
\overline{c}_{2,1}(A):K_{2}(A)\{\ell \}\rightarrow H_{cont}^{1}(A,\mathbb{Q}%
_{\ell }/\mathbb{Z}_{\ell }(2))/Div^{1}(A)  \label{cbar}
\end{equation}
where, for simplicity, $Div^{1}(A)$ denotes the maximal $\ell
-$divisible subgroup of $H_{cont}^{1}(A,\mathbb{Q}_{\ell
}/\mathbb{Z}_{\ell }(2))$. The exact sequence needed for this
factorization requires using continuous \'{e}tale cohomology which
Jannsen \cite{Jannsen} has developed. It agrees with continuous
group cochain cohomology when $A$ is a field and has all the needed
properties to construct the exact sequence in (\ref{c21}) and so to
show that $\overline{c}_{2,1}$ is properly defined. Moreover, for a
coefficient system of sheaves $F_{n}$ , when $H^{i-1}(A,F_{n})$
satisfies the Mittag-Leffler condition,
$H_{cont}^{i}(A,\lim\limits_{\leftarrow
}F_{n}) $ is the usual $\ell $-adic cohomology. Thus $H_{cont}^{1}(A,\mathbb{Z}%
_{\ell }(2))=H^{1}(A,\mathbb{Z}_{\ell }(2))$ and
$H_{cont}^{1}(A,\mathbb{Q}_{\ell }/\mathbb{Z}_{\ell
}(2))=H^{1}(A,\mathbb{Q}_{\ell }/\mathbb{Z}_{\ell }(2))$. In
addition, \cite[Theorem 5.14]{Jannsen} is the identification
$H^{2}_{cont}(A,\mathbb{Z}_{\ell }(2))\{\ell\}\cong
Im(H^{1}(A,\mathbb{Q}_{\ell}/\mathbb{Z}_{\ell }(2))).$

The proof of Theorem 3, while different in $char(A)>0$ and
$char(A)=0$, starts from the comparison of the Gersten-Quillen exact
sequence with a modification of the Bloch-Ogus sequence which is
still exact. Assume $A$ is a semi-local ring which is smooth and
essentially of finite type over a field $k$ ($k=\mathbb{Q}$ in
characteristic $0$). Consider, as in the proof earlier,
\begin{equation}
\begin{array}{ccc}
0\rightarrow K_{2}(A)\{\ell \}\rightarrow & K_{2}(F)\{\ell \}\rightarrow &
\coprod_{x\in Spec(A)^{1}}K_{1}(\kappa (x))\{\ell \} \\
\downarrow \overline{c}_{2,1}(A) & \downarrow \overline{c}_{2,1}(F) &
\downarrow = \\
0\rightarrow H^{1}(A,(2))/Div^{1}(A)\rightarrow & H^{1}(F,(2))/Div^{1}(F)%
\rightarrow & \coprod_{x\in Spec(A)^{1}}H^{0}(\kappa (x),(1))
\end{array}
\label{GQBO}
\end{equation}
where $H^{1}(A,(2))$ denotes $H^{1}(A,\mathbb{Q}_{\ell
}/\mathbb{Z}_{\ell }(2))$,
etc. for typographical reasons. If $char(A)$ is positive, then $H^{1}(A,\mathbb{%
Z}_{\ell }(2))$ is finite and so $H^{1}(A,\mathbb{Q}_{\ell }(2))$
which maps
onto the maximal divisible subgroup of $H^{1}(A,\mathbb{Q}_{\ell }/\mathbb{Z}%
_{\ell }(2))$ vanishes. Thus the bottom sequence in (\ref{GQBO}) is just the
Bloch-Ogus sequence and so $\overline{c}_{2,1}(A)$ is an isomorphism. If $%
char(A)=0$ and $A$ is essentially of finite type over $\mathbb{Q}$, then $%
H^{0}(\kappa (x),\mathbb{Q}_{\ell }/\mathbb{Z}_{\ell }(1))$ has no
divisible part for any $x\in Spec(A)^{(1)}$. Thus
$Div^{1}(A)=Div^{1}(F)$, and so the bottom sequence is an exact
quotient of the Bloch-Ogus sequence. Thus in either case
$\overline{c}_{2,1}(A)$ is an isomorphism.

Now suppose $A$ is an arbitrary semi-local ring of positive
characteristic. Since both $K_{2}(A)$ and $H^{1}(A,\mathbb{Q}_{\ell
}/\mathbb{Z}_{\ell }(2))$ commute with direct limits of rings we may
assume that $A$ is a semi-local ring essentially of finite type over
a field. We can then use Theorem \ref{fund} since
$H^{1}(A,\mathbb{Q}_{\ell }/\mathbb{Z}_{\ell }(2))$ is local for the
\'{e}tale topology with respect to closed sets by 1) in Example
(\ref{example}) and $K_{2}(A)\{\ell \}$ is locally of finite
presentation and satisfies the epic hensel pair condition. This
latter condition follows since $K_{3}(A;\mathbb{Z}/\ell ^{n})$ is
local for closed sets by 2) in Example (\ref{example}) and so
$_{\ell ^{n}}K_{2}(A),$ being a quotient, satisfies the epic hensel
pair condition. Note that the Corollary to Theorem \ref{fund} shows
that $K_{2}(A)\{\ell \}$ is then local for the \'{e}tale topology
with respect to closed sets.

The characteristic $0$ case is essentially the same but is
complicated by the fact that $H^{1}(A,\mathbb{Q}_{\ell
}/\mathbb{Z}_{\ell }(2))/Div^{1}(A)$ does not behave nicely with
respect to limits of rings. However Jannsen has shown that
\begin{equation*}
H^{1}(A,\mathbb{Q}_{\ell }/\mathbb{Z}_{\ell
}(2))/Div^{1}(A)\approxeq H_{cont}^{2}(A,\mathbb{Z}_{\ell
}(2))\left\{ \ell \right\} .
\end{equation*}
In addition he shows that the sequence,
\begin{equation*}
0\rightarrow \lim_{\leftarrow }{}^{1}H^{2}(B,\mathbb{Z}/\ell
^{n}(2))\rightarrow H_{cont}^{2}(B,\mathbb{Z}_{\ell }(2))\rightarrow
\lim_{\leftarrow }H^{2}(B,\mathbb{Z}/\ell ^{n}(2))\rightarrow 0,
\end{equation*}
is exact where $\lim_{\leftarrow }{}^{1}\,$ refers to the first left
derived functor of the system $\left\{ H^{2}(B,\mathbb{Z}/\ell
^{n}(2))\right\} .$ Consequently
\begin{equation*}
H_{cont}^{2}(B^{h},\mathbb{Z}_{\ell }(2))\approxeq
H_{cont}^{2}(A,\mathbb{Z}_{\ell }(2))
\end{equation*}
if $B^{h}$ is the henselization with respect to $I$ of a presentation $%
A\approxeq B/I$ as a quotient of a semi-local ring $B$ smooth over $k.$ This
and the corresponding isomorphism
\begin{equation*}
H^{1}(B^{h},\mathbb{Q}_{\ell }/\mathbb{Z}_{\ell }(2))\approxeq H^{1}(A,\mathbb{Q}%
_{\ell }/\mathbb{Z}_{\ell }(2))
\end{equation*}
then allow us to conclude that
\begin{equation*}
H_{cont}^{2}(B^{h},\mathbb{Z}_{\ell }(2))\left\{ \ell \right\}
\approxeq H_{cont}^{2}(A,\mathbb{Z}_{\ell }(2))\left\{ \ell \right\}
.
\end{equation*}
Now we may apply Theorem to show $\overline{c}_{2,1}:K_{2}(A)\{\ell
\}\rightarrow H_{cont}^{2}(A,\mathbb{Z}_{\ell }(2))\left\{ \ell
\right\} $ is an isomorphism.

Corollary 4 is obtained by noting that Suslin \cite[Corollary
2.8]{suslin} also proved that $H^{1}(F,\mathbb{Z}\mathbf{_{\ell
}}(2))$ is finite if $F$ is a field of positive characteristic. But
if $A$ is a unibranch local ring,
e.g. a normal local ring, then $H^{1}(A,\mathbb{Z}\mathbf{_{\ell }}%
(2))\hookrightarrow H^{1}(F,\mathbb{Z}\mathbf{_{\ell }}(2))$ is also
finite. Consequently $H^{1}(A,\mathbb{Q}\mathbf{_{\ell }}(2))=0$,
and so we have an
isomorphism $H^{1}(A,\mathbb{Q}\mathbf{_{\ell }/}\mathbb{Z}\mathbf{_{\ell }}%
(2))\cong H^{2}(A,\mathbb{Z}\mathbf{_{\ell }}(2))\{\ell \}$.

\appendix
\section{}

In this appendix $\zeta _{n}$ will always stand for a primitive
$n^{th}$ root of unity. We wish to prove the following theorem.
(Note that this is part of \cite[Theorem 1]{Kahn} and is proved
without assuming the generalized Kato conjecture.)

\begin{theorem}
Let $A$ be a connected, regular semi-local ring with quotient field $K,$ and
$n$ an odd positive integer relatively prime to the residue characteristics
of $A.$ Then
\begin{equation*}
H^{2}\left( A,\mathbb{Z}/n(2)\right) \overset{j}{\rightarrow }H^{2}\left( K,%
\mathbb{Z}/n(2)\right)
\end{equation*}
is injective.
\end{theorem}

We begin with a series of lemmas to calculate the cohomology of the
extension obtained by adjoining an $n^{th}$ root of $1$. Fix the following
situation in order to describe $\mathbb{Z}/n(k)$ and the action of the
cyclotomic character on this module. Fix an integer $n$ such that $n=dm$ and
$m=de$ for some positive integers $d,e,$ and let $A$ be a connected ring
with $\zeta _{m}\in A$ and $d$ a unit in $A.$ Assume $A\left[ \zeta
_{n}\right] $ is connected. Then
\begin{equation*}
X^{d}-\zeta _{m}=\prod_{i=1}^{d}\left( X-\zeta _{n}\left( \zeta
_{m}^{e}\right) ^{i}\right)
\end{equation*}
in $A\left[ \zeta _{n}\right] ,$ and $A\left[ \zeta _{n}\right] $ is a
Galois extension of $A$ with cyclic Galois group $G_{d}:=\mathbb{Z}%
/d=\left\langle \sigma \right\rangle \subset \left( \mathbb{Z}/n\right)
^{\times }.$

Construct a dictionary for the action of $G_{d}$ on $\mathbb{Z}/n(k)$ by $%
\sigma $ sends $\zeta _{n}\mapsto \zeta _{n}\left( \zeta _{m}^{e}\right) .$
Since $n\mid m^{2},$ we have
\begin{equation*}
\begin{tabular}{|l|l|}
\hline
$G_{d}$ on $\mathbb{Z}/n(k)$ & $\mathbb{Z}/d$ on $\mathbb{Z}/n$ \\ \hline
$\left( \zeta _{n},\zeta _{m}=\zeta _{n}^{d},\zeta _{d}=\zeta
_{m}^{e}\right) \in \mathbb{Z}/n(1)^{\times 3}$ & $\left( 1,d,m=de\right)
\in \mathbb{Z}/n^{\times 3}$ \\ \hline
$\sigma (\zeta _{n})=\zeta _{n}\left( \zeta _{m}^{e}\right) =\zeta
_{n}^{1+m}\in \mathbb{Z}/n(1)$ & $\sigma (1)=1+m\in \mathbb{Z}/n$ \\ \hline
$\sigma ^{k}(\zeta _{n})=\left( \zeta _{n}^{1+m}\right) ^{k}=\zeta
_{n}^{1+km}\in \mathbb{Z}/n(k)$ & $\sigma ^{k}(1)=\left( 1+m\right)
^{k}=1+km\in \mathbb{Z}/n$ \\ \hline
\end{tabular}
\end{equation*}
This dictionary makes the computation of the cohomology of $G_{d}$ on $%
\mathbb{Z}/n(k)$ straightforward.

\begin{lemma}
$H^{0}\left( G_{d},\mathbb{Z}/n(k)\right) =\mathbb{Z}/\left( m(d,k)\right) ,$
and

$H^{2r}\left( G_{d},\mathbb{Z}/n(k)\right) =\left\{
\begin{array}{lll}
\mathbb{Z}/(d,k) & {\text{ if}} & d{\text{ is odd or }}k\text{ is even or }%
4\mid m \\
\mathbb{Z}/2(d,k) & \text{if} & d\text{ is even and }k\text{ is odd and }%
4\nmid m
\end{array}
\right. $ and

$H^{2r+1}\left( G_{d},\mathbb{Z}/n(k)\right) =\left\{
\begin{array}{lll}
\mathbb{Z}/(d,k) & \text{if} & d\text{ is odd or }k\text{ is even or }4\mid m
\\
\mathbb{Z}/2(d,k) & \text{if} & d\text{ is even and }k\text{ is odd and }%
4\nmid m
\end{array}
\right. $
\end{lemma}

\begin{proof}
The standard resolution used to calculate the cohomology of a cyclic group $%
G_{d}$ on $\mathbb{Z}/n(k)$ is used to carry out this calculation. Our
dictionary allows us to calculate $N_{k}=\sum_{i=0}^{d-1}\sigma ^{ik}:%
\mathbb{Z}/n\rightarrow \mathbb{Z}/n$ and $T_{k}=1-\sigma ^{k}:\mathbb{Z}%
/n\rightarrow \mathbb{Z}/n$ as
\begin{eqnarray*}
T_{k}(j) &=&j-j\left( 1+km\right) =-kmj \\
N_{k}(j) &=&j\left( \sum_{i=0}^{d-1}\left( 1+ikm\right) \right) \\
&=&j\left( d+km\frac{d\left( d-1\right) }{2}\right) \\
&=&\left\{
\begin{array}{lll}
jd & \text{if} & d\text{ is odd or }k\text{ is even} \\
j\left( d-m\frac{d}{2}\right) & \text{if} & d\text{ is even and }k\text{ is
odd}
\end{array}
\text{ }\right.
\end{eqnarray*}
In the case when $d$ is even and $k$ is odd, we note that
\begin{equation*}
\left( m,1-\frac{m}{2}\right) =\left\{
\begin{array}{ll}
2 & \text{if }4\nmid m \\
1 & \text{otherwise}
\end{array}
\right. .
\end{equation*}
Then the sequence that calculates the cohomology groups is
\begin{equation*}
\cdots \mathbb{Z}/n\overset{N_{k}}{\rightarrow }\mathbb{Z}/n\overset{T_{k}}{%
\rightarrow }\mathbb{Z}/n\overset{N_{k}}{\rightarrow }\mathbb{Z}/n\overset{%
T_{k}}{\rightarrow }\mathbb{Z}/n\cdots ,
\end{equation*}
and we find
\begin{eqnarray*}
\text{Ker}\,T_{k} &=&\frac{d}{\left( d,k\right) }\mathbb{Z}/n\cong \mathbb{Z}%
/\left( m\left( d,k\right) \right) \\
\text{Im}\,T_{k} &=&km\mathbb{Z}/n=\left( d,k\right) \cdot m\mathbb{Z}/n \\
\text{Ker}\,N_{k} &=&\left\{
\begin{array}{lll}
m\mathbb{Z}/n & \text{if} & d\text{ is odd or }k\text{ is even or }4\mid m
\\
\frac{m}{2}\mathbb{Z}/n & \text{if} & d\text{ is even, }k\text{ is odd, and }%
4\nmid m
\end{array}
\text{ }\right. \\
\text{Im}\,N_{k} &=&\left\{
\begin{array}{lll}
d\mathbb{Z}/n & \text{if} & d\text{ is odd or }k\text{ is even or }4\mid m
\\
2d\mathbb{Z}/n & \text{if} & d\text{ is even, }k\text{ is odd, and }4\nmid m
\end{array}
\text{ }\right.
\end{eqnarray*}
The rest is straightforward.
\end{proof}

\begin{corollary}
\label{App:Cor} Let $A$ be a regular, connected semi-local domain, and let $%
n $ be an odd positive integer relatively prime to the residue
characteristics of $A$. Suppose that $m$ is the largest divisor of $n$ such
that $\zeta _{m}\in A$. Let $A_{n}=A[\zeta _{n}],$ and let $G=$Gal$(A_{n}/A).
$Then $H^{p}\left( G,\mathbb{Z}/n(2)\right) =\left\{
\begin{array}{lll}
\mathbb{Z}/m(2) & \text{if} & p=0 \\
0 & \text{if} & p>0
\end{array}
\right. $
\end{corollary}

\begin{proof}
If $m=n,$ we are done. Suppose that $n=\ell n^{\prime }$ where $\ell $ is
the smallest prime divisor of $\frac{n}{m}.$ Then there is a cyclic subgroup
$H\vartriangleleft G$ such that $A_{n^{\prime }}:=A^{H}=A[\zeta _{n^{\prime
}}],$ and $A_{n^{\prime }}$ is a Galois extension of $A$ with abelian Galois
group $G/H.$

If $\zeta _{\ell }\in A_{n^{\prime }},$ then $A_{n}$ and $A_{n^{\prime }}$
satisfy the hypothesis of the lemma. Hence
\begin{equation*}
H^{p}(H,\mathbb{Z}/n(2))=\left\{
\begin{array}{lll}
\mathbb{Z}/n^{\prime }(2) & \text{if} & p=0 \\
0 & \text{if} & p>0
\end{array}
\right. .
\end{equation*}
If $\zeta _{\ell }\neq A_{n^{\prime },}$ then $\left[ A_{n}:A_{n^{\prime
}}\right] =\ell -1$ is relatively prime to $n^{\prime }.$ Hence
\begin{equation*}
H^{p}(H,\mathbb{Z}/n(2)=\left\{
\begin{array}{lll}
\mathbb{Z}/n^{\prime }(2) & \text{if} & p=0 \\
0 & \text{if} & p>0
\end{array}
\right. .
\end{equation*}

Now by using the Hochschild-Serre spectral sequence,
\begin{equation*}
E_{2}^{p,q}=H^{p}\left( G/H,H^{q}\left( H,\mathbb{Z}/n\left( 2\right)
\right) \right) \Rightarrow H^{p+q}\left( G,\mathbb{Z}/n(2)\right) ,
\end{equation*}
and the above calculations, we conclude that $H^{p}\left( G/H,H^{0}\left( H,%
\mathbb{Z}/n\left( 2\right) \right) \right) \overset{\thickapprox }{%
\rightarrow }H^{p}\left( G,\mathbb{Z}/n\left( 2\right) \right) .$ Since $%
H^{0}\left( H,\mathbb{Z}/n\left( 2\right) \right) =\mathbb{Z}/n^{\prime
}(2), $ we may apply an induction hypothesis and thus prove the corollary.
\end{proof}

\begin{lemma}
\label{App:Lem} Let $A$ be a regular, connected, semi-local domain with
quotient field $K,$ $n$ a positive integer relatively prime to the residue
characteristics of $A$. Let $B$ be a Galois extension of $A$ with quotient
field $L$, and let $G=Gal(B/A).$ Suppose $\zeta _{n}\in B.$ Then
\begin{equation*}
H^{1}(G,H^{1}\left( B,\mathbb{Z}/n(2)\right) )\rightarrow
H^{1}(G,H^{1}\left( L,\mathbb{Z}/n(2)\right) )
\end{equation*}
is injective.
\end{lemma}

\begin{proof}
Since $A$ is a UFD, we have an exact sequence
\begin{equation*}
0\rightarrow A^{*}\rightarrow K^{*}\overset{\coprod ord_{x}}{\rightarrow }%
\coprod_{x\in A^{(1)}}i_{x*}\mathbb{Z}\rightarrow 0
\end{equation*}
where $A^{(1)}=\left\{ x\in Spec(A)/ht(x)=1\right\} $ and $%
ord_{x}:K^{*}\rightarrow i_{x*}\mathbb{Z}$ calculates the order of $f\in
K^{*}$ at the discrete valuation associated to $x,$ and $i_{x}:Spec(k(x))%
\rightarrow Spec(A)$ is the inclusion of the residue field of $A_{x}.$ If we
map this sequence to itself by multiplication by $n$, the cokernels form the
short exact sequence
\begin{equation*}
0\rightarrow A^{*}/A^{*n}\rightarrow K^{*}/K^{*n}\overset{\coprod ord_{x}}{%
\rightarrow }\coprod_{x\in A^{(1)}}i_{x*}\mathbb{Z}/n\rightarrow 0.
\end{equation*}
There is a sequence of $G$-modules obtained by twisting the sequence for $B$
with the cyclotomic character
\begin{equation*}
0\rightarrow B^{*}/B^{*n}(1)\rightarrow L^{*}/L^{*n}(1)\overset{}{%
\rightarrow }\coprod_{y\in B^{(1)}}i_{y*}\mathbb{Z}/n(1)\rightarrow 0,
\end{equation*}
relating $B/^{*}B^{*n}(1)=H^{1}\left( B,\mathbb{Z}/n(2)\right) $ to $%
L^{*}/L^{*n}(1)=H^{1}\left( L,\mathbb{Z}/n(2)\right) $ since $\zeta _{n}\in
L^{*}.$ Thus we need only show that
\begin{equation*}
H^{0}\left( G,L^{*}/L^{n*}(1)\right) \rightarrow H^{0}(G,\coprod_{y\in
B^{(1)}}i_{y*}\mathbb{Z}/n(1))
\end{equation*}
is onto to complete the proof of the lemma. But
\begin{equation*}
H^{0}(G,\coprod_{y\in B^{(1)}}i_{y*}\mathbb{Z}/n(1))\approxeq \coprod_{x\in
A^{(1)}}H^{0}\left( D_{y\mid x},i_{y*}\mathbb{Z}/n(1)\right) \approxeq
\coprod_{x\in A^{(1)}}i_{x*}\mathbb{Z}/m_{x}(1)
\end{equation*}
where $D_{y\mid x}$ is the decomposition group of a choice of $y\in B^{(1)}$
lying over $x\in A^{(1)}$ and $\mathbb{Z}/m_{x}(1)$ is the group of roots of
unity of order $n$ in $k(y)$ fixed by $D_{y\mid x}$ (expressed in the above
notation).

Any element $\alpha \in \coprod_{x\in A^{(1)}}i_{x*}\mathbb{Z}/m_{x}(1)$ can
be written as
\begin{equation*}
\alpha =\sum_{x\in A^{(1)}}\delta _{x}\cdot \zeta _{n}^{r_{x}}
\end{equation*}
where the sum is finite and $\delta _{x}=0$ for $x^{\prime }\in A^{(1)}$ if $%
x^{\prime }\neq x$ and $\delta _{x}=1\in i_{x*}\mathbb{Z}/n$ at $x\in
A^{(1)}.$ Then there is $f_{x}\in K^{*}/K^{*n}$ with $ord_{x}(f_{x})=\delta
_{x}$ and $ord_{x^{\prime }}(f_{x})=0$ if $x^{\prime }\neq x.$ Hence $%
\sum_{x\in A^{(1)}}\left( \sum_{\sigma \in G/D_{y\mid x}}f_{x}\sigma \left(
\zeta _{n}^{r_{x}}\right) \right) $ has image $\alpha .$
\end{proof}

\begin{remark}
The proof amounts to first reducing to the case of a dvr $A$ by using the
UFD sequence and then solving the problem at the completion where $D_{y\mid
x}$ becomes the Galois group and finally spreading that solution around
using the transitivity of the action of $G$ on $\left\{ y\in B^{(1)}/y\text{
lies over a fixed }x\in A^{(1)}\right\} .$
\end{remark}

\begin{proof}
If $\zeta _{n}\in A,$ the theorem reduces to the well known result $%
Br(A)\subset Br(K)$ since $H^{2}\left( A,\mathbb{Z}/n(k)\right) \approxeq
H^{2}\left( A,\mathbb{Z}/n(1)\right) \approxeq \,_{n}Br(A)$ and similarly
for $K.$ We reduce the theorem to this case by analyzing the cohomology of a
cyclic Galois covering $A\left[ \zeta _{n}\right] /A$ with Galois group $G.$
There is a spectral sequence for this covering
\begin{equation*}
E_{2}^{p,q}=H^{p}\left( G,H^{q}\left( A[\zeta _{n}],\mathbb{Z}/n\left(
2\right) \right) \right) \Rightarrow H^{n}\left( A,\mathbb{Z}/n(2)\right) .
\end{equation*}
Suppose $x\in Ker\left[ H^{2}(A,\mathbb{Z}/n(2))\rightarrow H^{2}(K,%
\mathbb{Z}/n(2))\right] .$ Then $x=0$ in $H^{2}\left( A[\zeta _{n}],%
\mathbb{Z}/n\left( 2\right) \right) ^{G}$ by the above remarks. But Lemma
\ref{App:Lem} shows that $x=0$ also in $E_{\infty }^{1,1}\subseteq
H^{1}\left( G,H^{1}\left( A[\zeta _{n}],\mathbb{Z}/n\left( 2\right) \right)
\right) .$ Finally Corollary \ref{App:Cor} shows that $x\in E_{\infty
}^{2,0}\subseteq H^{2}\left( G,\mathbb{Z}/m\left( 2\right) \right) $ where $%
\zeta _{m}\in A.$ But $H^{2}\left( G,\mathbb{Z}/m\left( 2\right) \right) $
is the same for the Galois covering $A[\zeta _{n}]/A$ and $K[\zeta _{n}]/K. $
Hence $x=0$ and the theorem follows.
\end{proof}


\begin{thebibliography}{99}
\bibitem{blochogus}  S. Bloch and A. Ogus ``Gersten's conjecture and the
homology of schemes'', Ann Sci Ecole Norm Sup., 4th s\'{e}r. 7(1974).
181-202.

\bibitem{CTO}  J-L Colliot-Th\'{e}l\`{e}ne and M. Ojanguren ``Espaces
Principaux homog\`{e}nes localement triviaux'', Publications Mathematiques
de l'IHES, 75(1992), 97-122.

\bibitem{CHK}  J.-L. Colliot-Th\'{e}l\`{e}ne, R.T. Hoobler, and B. Kahn,
``The Bloch---Ogus---Gabber Theorem'', \underline{Algebraic} \underline{%
K-Theory}, Fields Inst. Commun., Vol. 16, Amer. Math. Soc., Providence, RI,
1997, 31--94.

\bibitem{gabber}  O. Gabber, ``K-theory of henselian local rings and
henselian pairs'', a letter from Ofer Gabber to M. Karoubi, 83, AMS Contemp
Math, 1989.

\bibitem{gabberhensel}  O. Gabber, ``Affine analog of the proper base change
theorem'', Israel J. of Math., 87(1994),325-335.

\bibitem{gillet}  H. Gillet, ``Gersten's conjecture for the K-theory with
torsion coefficients of a discrete valuation ring'', J Algebra, 103(1986),
377-380.

\bibitem{gilletk}  H. Gillet,.and M. Levine, ``The relative form of
Gersten's conjecture over a discrete valuation ring:The smooth case'', J
Pure Appl Algebra, 46(1987), 59-71.

\bibitem{gilletbo}  H. Gillet, personal communication, July 1995.

\bibitem{Grayson}  D. Grayson, ``Universal exactness in algebraic
K-theory'', J Pure Appl Algebra, 36(1985), 139-141.

\bibitem{Jannsen}  U. Jannsen, ``Continuous \'{e}tale cohomology'', Math
Ann, 280(1988), 207-245.

\bibitem{Kahn}  B. Kahn, ``Deux th\'{e}or\`{e}mes de comparision en
cohomologie \'{e}tale; applications'', Duke Math J, 69(1993), 137-165.

\bibitem{Milne}  J. Milne, \underline{\'{E}tale} \underline{Cohomology},
Princeton University Press, 1980, Princeton, NJ.

\bibitem{Weibel}  C. Pedrini, and C. Weibel, ``Invariants of real curves'',
Rend Sem Mat Univers Politecn Torino, 49(1991), 139-173.

\bibitem{Quillen}  D. Quillen, ``Higher Algebraic K-theory I'', \underline{%
Algebraic} \underline{K-theory I}: \underline{Higher} \underline{K-theories}%
, Lecture Notes in Math, 341, Springer-Verlag, Berlin, 1973, 85-147.

\bibitem{Raynaud}  Michel Raynaud, \underline{Anneaux} \underline{Locaux}
\underline{Hens\'{e}lien} , Lecture Notes in Math, 169, Springer-Verlag,
Berlin, 1970.

\bibitem{sga4}  \underline{SGA4}:\underline{Theorie} \underline{des}
\underline{Topos} \underline{et} \underline{Cohomologie} \underline{\'{E}tale%
} \underline{des} \underline{Schemas}, Lecture Notes in Math, 305,
Springer-Verlag, Berlin, 1973.

\bibitem{soule}  C. Soul\'{e}, ``K-th\'{e}orie des anneaux d'entiers de
corps de nombres et cohomologie \'{e}tale'', Inventiones math, 55(1979),
251-295.

\bibitem{Strano}  R. Strano,``On the \'{e}tale cohomology of hensel rings'',
Comm in Alg, 12, 2195-2211.

\bibitem{suslin}  A. Suslin, ``Torsion in $K_{2}$ of fields'', K-theory,
1(1987), 5-29.

\bibitem{Tate}  J. Tate, ``Relations between $K_{2}$ and Galois
cohomology'', Inv Math 36(1976), 257-274.

\bibitem{LLouise}  C. Weibel, ``\'{E}tale Chern classes at the prime 2'',
\underline{Proceedings} \underline{of} \underline{the} \underline{Lake}
\underline{Louise} \underline{Conference} \underline{on} \underline{Algebraic%
} \underline{K-Theory}, NATO ASI Series C, 407, Kluwer Academic Publishers,
Dordrecht 1993, 249-286.
\end{thebibliography}
\end{document}